\documentclass[12pt,reqno]{amsart}
\usepackage{latexsym,amsmath}
\usepackage{amsthm}
\usepackage{amssymb}
\usepackage{amsthm}
\usepackage{epsfig}
\begin{document}
\newtheorem{theorem}{Theorem}[section]
\newtheorem{cor}[theorem]{Corollary}
\newtheorem{lemma}[theorem]{Lemma}
\newtheorem{fact}[theorem]{Fact}
\newtheorem{claim}[theorem]{Claim}
\newtheorem{definition}[theorem]{Definition}
\theoremstyle{definition}
\newtheorem{example}[theorem]{Example}
\newtheorem{remark}[theorem]{Remark}
\newcommand\eps{\varepsilon}
\newcommand{\pp}{\mathcal P}
\newcommand{\E}{\mathbb E}
\newcommand{\Var}{{\rm Var}}
\newcommand{\Prob}{\mathbb{P}}
\newcommand{\w}{\mathbf{w}}
\newcommand{\N}{{\mathbb N}}
\newcommand{\R}{{\mathbb R}}
\newcommand{\eqn}[1]{(\ref{#1})}

\newcommand\qedc{\hfill $\Diamond$ \smallskip}
\renewcommand{\labelenumi}{(\roman{enumi})}

\def\bibitemx#1{\bibitem{#1}\marginpar{$\otimes$}}
\def\co{~\marginpar{\fbox{,}}}
\def\wck#1 {\underline{#1}~\marginpar{\fbox{#1} {\tiny ?}}}
\def\silent#1\par{\par}
\def\text#1{\quad\mbox{#1}\quad}

\def\quest#1{\footnote{$\rightarrow$ {\sf #1}}}

\makeatletter
\renewcommand{\@seccntformat}[1]{\@nameuse{the#1}.\quad}
\makeatother

\def\wepn{{\em wep}$_n$}

\title{Rank-based attachment leads to power law graphs}

\author{Jeannette Janssen and Pawe{\l} Pra{\l}at}
\thanks{The authors are supported by MITACS and NSERC. This work is part of the MITACS project Modelling and Mining of Networked Information Spaces (MoMiNIS)}

\address{Department of Mathematics and Statistics, Dalhousie University, Halifax NS, Canada B3H 3J5}
\email{\tt janssen@mathstat.dal.ca}

\address{Department of Mathematics and Statistics, Dalhousie University, Halifax NS, Canada B3H 3J5}
\email{\tt pralat@mathstat.dal.ca}


\begin{abstract}
We investigate the degree distribution resulting from graph generation models based on rank-based attachment. In rank-based attachment, all vertices are ranked according to a ranking scheme. The link probability of a given vertex is proportional to its rank raised to the power $-\alpha$, for some $\alpha\in(0,1)$. Through a rigorous analysis, we show that rank-based attachment models lead to graphs with a power law degree distribution with exponent $1+1/\alpha$ whenever vertices are ranked according to their degree, their age, or a randomly chosen fitness value. We also investigate the case where  the ranking is based on the initial rank of each vertex; the rank of existing vertices only changes to accommodate the new vertex. Here, we obtain a sharp threshold for power law behaviour. Only if initial ranks are biased towards lower ranks, or chosen uniformly at random, we obtain a power law degree distribution with exponent $1+1/\alpha$. This indicates that the power law degree distribution often observed in nature can be explained by a rank-based attachment scheme, based on a ranking scheme that can be derived from a number of different factors; the exponent of the power law can be seen as a measure of the strength of the attachment.
\end{abstract}

\maketitle


\section{Introduction}

The occurrence of power law degree distributions in self-organizing networks such as the web graph is often explained by a model based on the principle of \emph{Preferential Attachment (PA)}. In the original PA model proposed by Barab\'asi and Albert~\cite{BA}, new vertices join a graph one by one, and each new vertex chooses a pre-determined number of neighbours at random, so that the probability that a vertex is chosen as a neighbour (its {\sl link probability}) is proportional to its degree. Analysis shows that this model indeed generates power law graphs with high probability, where the exponent of the power law equals 3~\cite{BA,BRST}. More general PA models, such as the ones proposed and analyzed in~\cite{ACL} and~\cite{CF} allow for the creation of edges between existing vertices and the deletion of edges and vertices. The power law of the degree distribution in this case depends on the probabilities with which various kind of steps (edge addition, vertex addition, deletion) are taken. By varying these probabilities, any exponent in the range $(2,\infty)$ can be obtained. In~\cite{BCDR}, the \emph{preferential attachment with fitness} is studied, in which the degree of a vertex is scaled by its individual fitness factor to determine its attractiveness. 

It is tempting to conjecture that the exponent of the power law can be controlled by varying the strength of preferential attachment. Precisely, the link probability could be proportional to the degree raised to a power $\alpha$; the hope is that the exponent of the power law would be a function of $\alpha$. This would give a more natural way to fit the model to real-life data than that given by the generalized addition/deletion models: the ratio of edge vs.~vertex addition steps may well be dictated by considerations about the data that are independent of the power law. Unfortunately, as pointed out in~\cite{KR}, this approach does not work. Only the case where $\alpha =1$ (i.e.~the standard PA model) leads to a power law degree distribution. If $\alpha<1$ (weak preferential attachment), the degree distribution is a stretched exponential, while if $\alpha >1$ (strong preferential attachment), the graph will be close to a star, with one vertex adjacent to almost all other vertices.

In this paper, we show that the approach outlined above does give the proposed results if the preferential attachment is based on a \emph{ranking} of the vertices. In other words, the vertices are ranked from 1 to $n$ according to their degree (so the vertex with highest degree has rank 1, etc.), and the link probability of a given vertex is proportional to its rank, raised to the power $-\alpha$ for some $\alpha\in(0,1)$; we will refer to $\alpha$ as the \emph{attachment strength}. (Negative powers are chosen since a low value for rank should result in a higher link probability.) Then, with high probability, the resulting graph will have a power law with exponent $1+1/\alpha$. The rank-based approach was first proposed by Fortunato, Flammini and Menczer in~\cite{FFM}, and the occurrence of a power law was postulated based on simulations.

In~\cite{FFM}, the attachment strength $\alpha$ is allowed to be any positive real number. However, if $\alpha >1$, then only a vanishingly small proportion of all vertices have any acquired links at all. This is easy to show for the scenario where vertices are ranked according to age. We feel confident that the same holds for the other ranking schemes that lead to power law degree distributions. Since a scenario where almost all vertices have the minimum degree does not correspond to a typical self-organizing network, we have limited our analysis to the case where $\alpha\in(0,1)$.

As we will show, rank-based attachment leads to power law graphs for a variety of different ranking schemes. One obvious ranking scheme is to rank vertices by age (\emph{the old get richer}); we show that this leads to a power law with the same exponent $1+1/\alpha$. A more general graph model with rank-based attachment based on ranking according to age is the \emph{protean graph} model, which was proposed and explored by \L{}uczak, Pra\l{}at and Wormald in~\cite{TLPP, PPNW, PP}. It is also discussed in~\cite{FFM} and~\cite{GKK}. As a contrast, in this paper we also consider an inverse age ranking scheme, where younger vertices are ranked higher. As can be expected, this scheme is not likely to lead to a heavy tail degree distribution: with high probability, the maximum degree is of order $\log n$, where $n$ is the total number of vertices.

In~\cite{FFM}, a ranking scheme based on an external prestige label for each vertex is given, and it is shown through a heuristic analysis and simulations that this scheme also leads to power law graphs, with the same exponent. Precisely, each vertex at its birth is assigned a randomly chosen {\sl fitness} $\ell\in (0,1)$, and vertices are ranked according to their prestige label. As argued below, since the ranking is based only on the relative values of the fitness values, the distribution according to which $\ell$ is chosen is irrelevant (all distributions give equivalent graph generation processes).

In order to allow for a different distribution of \lq\lq prestige\rq\rq  over the vertices, we considered also a \emph{random ranking} scheme. Here, each vertex is assigned an initial rank according to a given distribution. We consider distributions of the following form. Let $R_i$ be the initial rank of a vertex born at time $i$. Then $\Prob (R_i\leq k)= (k/i)^s$.  First we show that, if $s=1$, then the situation is similar to the one described previously, and vertices with initial rank $R_i$ exhibit behaviour as if they had received fitness $R_i/i$. Thus, we obtain a power law graph.

Next, we consider the case where $s>1$, so the rank of new vertices is biased towards the lower ranks (note that low rank refers to a vertex with high value of $R_i$). In this case, with high probability the rank of a vertex will remain close to its original value throughout the process, so the behaviour is similar to the case of ranking according to age, and we obtain a power law graph. If $s<1$, so initial ranks are biased towards the higher ranks, then we show that vertices tend to drift towards the lower ranks, and the behaviour is similar to that of ranking according to inverse age, where no power law is likely to occur. Thus, the value of $s=1$ gives a sharp threshold for power law behaviour of the degree distribution.

These results suggest an explanation for the power law degree distribution often observed in real-life networks such as the web graph, protein interaction networks, and social networks. The growth of such networks can be seen as governed by a rank-based attachment scheme, based on a ranking scheme that can be derived from a number of different factors such as age, degree, or fitness. The exponent of the power law is independent of these factors, but is rather a consequence of the attachment strength. In addition, rank-based attachment accentuates the difference between higher ranked vertices: the difference in link probability between the vertices ranked 1 and 2 is much larger than that between the vertices ranked 100 and 101. This again corresponds to our intuition of what constitutes a credible mechanism for link attachment.

In order to establish the right attachment strength to model a given real-life network we should consider the following. In a graph in which the number of vertices of degree $k$ decreases roughly as $k^{-\gamma}$ the fraction of vertices of degree at least $k$ changes roughly as
$$
\sum_{\ell\ge k}O(\ell^{-\gamma})=O(k^{1-\gamma})\,.
$$
Thus, in order to imitate this distribution the attachment strength $\alpha$ should be set to  $\alpha\sim 1/(\gamma-1)$. For the web graph the distribution of total degrees is, at this moment, unknown. However, the number of vertices of in-degree $k$ decreases roughly as $k^{-2.1}$, while the fraction of vertices of out-degree $k$ can be approximated by $k^{-2.7}$ (cf., Broder {\em et al.}~\cite{Broder}). Thus, if the total degree of the graph is close to the distribution of in-degree, then the appropriate attachment strength for a rank-based model is $\alpha \sim 0.91$.


\section{Definitions}

In this section, we formally define the graph generation model based on rank-based attachment. The model produces a sequence $\{G_t\}_{t=1}^\infty = \{(V_t,E_t)\}_{t=1}^\infty$ of undirected graphs, where $t$ denotes time. Our model has two fixed parameters: initial degree $d\in N$, and attachment strength $\alpha \in (0,1)$. At each time $t$, each vertex $v\in V_t$ has  rank $r(v,t) \in [t]$ (we use $[t]$ to denote the set $\{1,2,\dots, t\}$). To obtain a proper ranking, the rank function $r(\cdot,t):V_t\rightarrow [t]$ is a bijection for all $t$, so every vertex has a unique rank. In agreement with the common use of the word \lq\lq rank\rq\rq, high rank refers to a vertex $v$ for which $r(v,t)$ is small: the highest ranked vertex is ranked number one, so has rank equal to 1; the lowest ranked vertex in $V_n$ has rank $n$. The initialization and update of the ranking is done according to a \emph{ranking scheme}. Various ranking schemes are considered in this paper; we first give the general model, and then list the ranking schemes.

For any $0 < \alpha < 1$, we define the function $g_\alpha:\N\rightarrow \R$:
\begin{equation}\label{denominator}
g_\alpha(t)=\sum_{j=1}^t j^{-\alpha}=\frac {t^{1-\alpha}}{1-\alpha} + O(1) \,.
\end{equation}

Let $G_1=(V_1,E_1)=(\{v_1\},\emptyset)$ be a fixed initial graph with a single vertex with $d$ loops, and rank $r(v_1,1)=1$. For $t>1$ we form $G_t$ from $G_{t-1}$ according to the following rules:
\begin{itemize}
     \item Add a new vertex $v_t$ together with $d$ edges from $v_t$ to existing vertices chosen randomly with weighted probabilities. The edges are added in $d$ substeps. In each substep, one edge is added, and the probability that $v_i$ is chosen as its endpoint (the link probability), equals $r(v_i,t-1)^{-\alpha} / g_\alpha(t-1)$.
     \item Update the ranking function $r(\cdot, t) : V_t \to [t]$ according to the ranking scheme.
 \end{itemize}

Our model allows for loops and multiple edges; there seems no reason to exclude them. However, there will not in general be very many of these, so excluding them can be shown not to affect our conclusions in any significant way.

We now define the different ranking schemes that are considered in this paper.

\begin{itemize}
   \item {\bf Ranking by age}:
The vertex added at time $t$ obtains a rank $t$ and retains this rank. That is, $r(v_i,t) = i$ for $i \in [t]$.
   \item {\bf Ranking by inverse age}:
The vertex added at time $t$ obtains a rank $1$, but its rank shifts by one each time a new vertex is added. That is, $r(v_i,t) = t-i+1$ for $i \in [t]$.
   \item {\bf Ranking by random labeling}:
The vertex added at time $t$ obtains a label $l(v_t)\in (0,1)$ chosen uniformly at random. Vertices are ranked according to their labels: if $l(v_i)<l(v_j)$, then $r(v_i,t)<r(v_j,t)$.
   \item {\bf Random ranking}:
The vertex added at time $t$ obtains an initial rank $R_t$ which is randomly chosen from $[t]$ according to a prescribed distribution. Formally, let $F:[0,1]\rightarrow [0,1]$ be any cumulative distribution function. Then for all $k\in[t]$,
\[
\Prob (R_t\leq k)=F(k/t).
\]
   \item {\bf Ranking by degree}:
After each time step $t$, vertices are ranked according to their degrees in $G_t$, and ties are broken by age. Precisely, if $\deg(v_i,t)>\deg(v_j,t)$ then $r(v_i,t)<r(v_j,t)$, and if $\deg(v_i,t)=\deg(v_j,t)$ then $r(v_i,t)<r(v_j,t)$ if $i<j$.
\end{itemize}

In the rest of the paper, $\{ G_t \}_{t=1}^{\infty}$ is assumed to be a graph sequence generated by the rank-based attachment model, with ranking scheme as defined in each particular section, and $d$ and $\alpha$ are assumed to be the initial degree and attachment strength parameters of the model as defined above. The results are generally about the degree distribution in $G_n$, where the asymptotics are based on $n$ tending to infinity.

We will use the stronger notion of \emph{wep} in favour of the more commonly used \emph{aas}, since it simplifies some of our proofs. We say that an event holds \emph{with extreme probability} (\emph{wep}), if it holds with probability at least $1-\exp (-\Theta(\log ^2 n))$ as $n \to \infty$.  Thus, if we consider a polynomial number of events that each holds \emph{wep}, then \emph{wep} all events hold. To combine this notion with asymptotic notations such as $O()$ and $o()$, we follow the conventions in~\cite{GTH}.


\section{Ranking by degree}

The first ranking scheme we consider is the ``ranking'' version of preferential attachment: vertices with higher degree are ranked higher. That is, the rank function $r(\cdot,t):V_t \to [t]$ is determined by the degree sequence at time $t$: if $\deg(v_i,t)>\deg(v_j,t)$, then $r(v_i,t)<r(v_j,t)$; otherwise (that is, if $\deg(v_i,t)=\deg(v_j,t)$) $r(v_i,t)<r(v_j,t)$ if $i<j$. In contrast to all other ranking schemes, this means that the rank of a vertex can change by more than one in each step; if the degree of a vertex changes from $k$ to $k+1$, then the change in its rank can be as large as the total number of vertices of degree $k$ and $k+1$.

In this section only, in order to omit tedious details, we assume that $d=1$. The general case can be studied in a similar way.

For all $t\geq 1$ and $k \ge 1$, let $Y_k(t)$ denote the number of vertices of degree $k$ in $G_t$, and let $Y_{\le k}(t)=Y_1(t)+\dots+Y_k(t)$. (Note that $Y_0(t) = 0$ for all $t \ge 1$.) Note that, at time $t$, the vertices of degree $k$ have ranks starting at $t-Y_{\le k}(t)+1$, and ending at $t-Y_{\le k-1}(t)$. When the edge $v_j v_{t+1}$ is added at time $t+1$, the change in the $Y_i$'s has contributions from two sources, namely, the change in degree of vertex $v_j$, and the addition of vertex $v_{t+1}$ of degree $1$. Note that the probability that a vertex of degree $k$ receives a link in step $t+1$ equals
\[
\sum_{j=t-Y_{\le k}(t)+1}^{t-Y_{\le k-1}(t)}\frac{j^{-\alpha}}{g_\alpha(t)}=\frac{g_\alpha(t-Y_{\le k-1}(t))-g_\alpha(t - Y_{\le k}(t))}{g_\alpha(t)}.
\]
Thus, the following equations express the expected change in each time step:
\begin{eqnarray*}
    \E (Y_1(t+1) - Y_1(t) ~|~ G_t) &=& 1 - \frac {g_\alpha(t)- g_\alpha(t-Y_1(t))} {g_\alpha(t)}
\end{eqnarray*}
and similarly, for all $k\ge 2$,
\begin{eqnarray*}
    \E (Y_k(t+1) - Y_k(t) ~|~ G_t) &=& \frac {g_\alpha(t - Y_{\le k-2}(t)) - g_\alpha(t - Y_{\le k-1}(t))}{g_\alpha(t)}\\
&&\quad - \frac {g_\alpha(t - Y_{\le k-1}(t))-g_\alpha(t - Y_{\le k}(t))}{g_\alpha(t)}\,.
\end{eqnarray*}
(Note that $Y_k(t)=0$ for all $k>t$.)

To analyze the behaviour of the $Y_i$, we use the \emph{differential equations method} (see for a general survey~\cite{NW}). First, by interpolating variables $Y_i(t)$ by real functions and \emph{presuming} that the changes in the functions are equal to their expected changes, the equations above can be turned into differential equations. The nature of the limiting behaviour as $n \to \infty$ can be emphasised by considering real functions $z_i(x)$ to model the behaviour of $\frac 1n Y_i(xn)$. Using the approximation~(\ref{denominator}), we obtain a system of differential equations:
\begin{eqnarray*}
z_1'(x) &=& \left( 1 - \frac {z_1(x)}{x} \right)^{1-\alpha} \\
z_k'(x) &=& \left( 1 - \frac {\sum_{j=1}^{k-2} z_j(x)}{x} \right)^{1-\alpha} -2 \left( 1 - \frac {\sum_{j=1}^{k-1} z_j(x)}{x} \right)^{1-\alpha}\\
&& + \left( 1 - \frac {\sum_{j=1}^{k} z_j(x)}{x} \right)^{1-\alpha} \mbox{\ \ \ \  for }k\ge 2\,,
\end{eqnarray*}
where $z_0(x) = 0$ for all $x$. The initial conditions are $z_k(0)=0$ for all $k$ (or more precisely, the right limit as $x$ approaches zero equals zero). The solutions of these equations are $z_k(x) = c_k x$, where the constants $c_k$ are defined below. Let $C_k = 1-\sum_{j=1}^{k} c_j=\sum_{j=k+1}^{\infty} c_j$, so $c_k=C_{k-1}-C_{k}$. Solving the differential equations one by one, we get the following recurrence:
\begin{eqnarray}
C_0 &=& 1 \nonumber \\
C_1 &=& 1 - C_1^{1-\alpha} \nonumber \\
C_{k-1} - C_{k} &=& C_{k-2}^{1-\alpha} - 2 C_{k-1}^{1-\alpha} + C_k^{1-\alpha}\mbox{\ \ \ \  for }k\ge 2\,. \label{eqn:c_k}
\end{eqnarray}
The recurrence~(\ref{eqn:c_k}) is telescoping, so it can be simplified.  Indeed,
\begin{eqnarray*}
C_k(1+C_k^{-\alpha}) &=& C_{k-1}^{1-\alpha} - C_{k-2}^{1-\alpha} + C_{k-1}(1+C_{k-1}^{-\alpha}) \\
&=& C_{k-1}^{1-\alpha} - C_{1}^{1-\alpha} + C_{2}(1+C_{2}^{-\alpha}) \\
&=& C_{k-1}^{1-\alpha}\,.
\end{eqnarray*}
Thus, $C_0 = 1$ and $C_k + C_k^{1-\alpha} = C_{k-1}^{1-\alpha }$ for $k \ge 1$.

\begin{lemma}\label{lemma:c_k}
Let the sequence $(C_k:k\geq 0)$ be recursively defined so that $C_0=1$, and $C_k$ is the unique positive solution to the equation
\begin{equation}
\label{form:C_k}
C_k + C_k^{1-\alpha} = C_{k-1}^{1-\alpha }.
\end{equation}
Then
$C_k= c_\alpha k^{-1/\alpha}(1+o(1))$, where $c_\alpha=\left(\frac{1-\alpha}{\alpha}\right)^{1/\alpha}$.
\end{lemma}

\begin{proof}
Note first that the function $f(x)=x+x^{1-\alpha}$ is concave and strictly increasing when $x>0$,  so $C_k$ is well-defined and decreasing. Let $B_k=C_k k^{1/\alpha}$. We will show that $B_k$ converges to a constant when $k\rightarrow\infty$. Substituting the expression for $B_k$ into (\ref{form:C_k}) we obtain the following recurrence relation:
\[
B_k^{1-\alpha}+\frac{B_k}{k}=\left(1+\frac{1}{k-1}\right)^{\frac{1-\alpha}{\alpha}} B_{k-1}^{1-\alpha}.
\]
Consider the function $f_k: [0,\infty)\rightarrow\R$ given by
\[
f_k(x)=x^{1-\alpha}+\frac{x}{k}-\left( 1+\frac{1}{k-1}\right)^{\frac{1-\alpha}{\alpha}}x^{1-\alpha}.
\]
The roots of the function $f_k$ are $x=0$ and $x=x_k$, where
\[
x_k = \left( k\left( 1+\frac{1}{k-1}\right)^{\frac{1-\alpha}{\alpha}}-k\right)^{\frac{1}{\alpha}}.
\]
Using the Taylor expansion of $\left( 1+\frac{1}{k-1}\right)^{\frac{1-\alpha}{\alpha}}$, and considering $x_k^\alpha$, we can show that $x_k$ is a decreasing sequence for $k\geq 2$, with limit $\left(\frac{1-\alpha}{\alpha}\right)^{1/\alpha}=c_{\alpha}$.

It is straightforward to verify that $f_k'(x_k)>0$. Since $f_k(x_k)=0$ and $x^{1-\alpha}+x/k$ is an increasing function of $x$, it follows that $B_{k-1}>x_k$ implies that $B_k<B_{k-1}$, and $B_{k-1}<x_k$ implies that $B_k>B_{k-1}$. Using the recursive expression for $B_k$, we obtain that
\begin{eqnarray}
\label{fkBk}
f_k(B_k)&=&B_k^{1-\alpha}+\frac{B_k}{k}-\left( 1+\frac{1}{k-1}\right)^{\frac{1-\alpha}{\alpha}}B_k^{1-\alpha}\notag\\
&=& \left( 1+\frac{1}{k-1}\right)^{\frac{1-\alpha}{\alpha}}(B_{k-1}^{1-\alpha}-B_k^{1-\alpha})
\end{eqnarray}

Since $f_k$ is increasing, $x>x_k$ if and only if $f_k(x)>f_k(x_k)$. Thus, if $B_{k-1}>x_k$ then $B_{k-1}>B_k$, so $f_k(B_k)>0$ and thus $B_k>x_k>x_{k+1}>c_\alpha$. Therefore, if there exists a value $k_0>2$ so that $B_{k_0-1}>x_{k_0}$, then $\{ B_k:k\geq k_0\}$ is monotone decreasing, and bounded from below by $c_\alpha$. If no such $k_0$ exists, then $B_{k-1}\leq x_k$ for all $k\geq 2$, and thus $B_k$ is monotonically increasing, and bounded from above by a converging sequence. In both cases, the sequence $B_k$ converges.

From~(\ref{fkBk}), we can then conclude that $f_k(B_k)\rightarrow 0$, and thus $B_k\rightarrow x_k\rightarrow c_{\alpha}$ as $k\rightarrow \infty$
\end{proof}

Since $c_k=C_{k-1}-C_k$, the above lemma implies that $c_k=(1/\alpha)c_\alpha k^{-(1+1/\alpha)}$.
We will see in the rest of this section that the variables $Y_i$ indeed follow the behaviour suggested by the discussion above, as expressed in the following theorem.

Let
\begin{equation*}
K = K(n) = \frac {n^{\alpha/(4\alpha+2)} } {\log^{2 \alpha / (2\alpha+1)} n}.
\end{equation*}

\begin{theorem}
\label{thm:degree}
\emph{Wep} for all $k, 1 \le k \le K$ and $t, 1\le t\le n$
$$
Y_k(t) =  c_k n (1+o(1)) = (1/\alpha)c_{\alpha} k^{1/\alpha +1} n(1+o(1)).
$$
\end{theorem}

Note that if the random variable $Y_k(t)$ is smaller than $c_k t$ at some point of the process, then the probability that $Y_k(t)$ increases in the next step will go up, so the value will be increasing with a higher rate. Likewise, if $Y_k > c_k t$, then the probability that $Y_k(t)$ stays equal goes up, so $Y_k(t)$ tends to stay at the same value for longer. Thus, these random variables have a certain ``self-repairing'' quality, so we expect them to behave well.

This suggests that, in order to show a concentration of $Y_k$, the differential equations method can be used. In this case, the full force of this method need not be used, but it suffices to use martingales, or bound the variables by stochastically dominating the behaviour using binomially distributed variables --- these can be considered as primitive versions of the differential equations method. We present a technique based on a well-known Azuma-Hoeffding inequality (see for example Lemma~4.2 in~\cite{NW}).

\begin{lemma}\label{lem:martingales}
Let $X_0, X_1, \dots, X_t$ be a supermartingale such that $|X_j-X_{j-1}| \le c_j$, $1 \le j \le t$, for constants $c_j$, Then for any $\xi>0$
$$
\Prob(X_t-X_0 \ge \xi) \le \exp \left(- \ \frac {\xi^2}{2 \sum c_j^2} \right)\,.
$$
\end{lemma}

To avoid tedious repetition, we present the full proof of Theorem~\ref{thm:degree} for variable $Y_1$ only. A proof sketch will be given for the other cases.

\begin{theorem}
\emph{Wep} for all $t, 1\le t\le n$
$$
Y_1(t) = c_1 t + O(\sqrt{n} \log n)\,.
$$
\end{theorem}
\begin{proof}
Let $X(t)=Y_1(t)-c_1t$. Since we expect $Y_1(t)$ to stay around $c_1 t$, and thus $X(t)$ to be a random variable close to a martingale. For any two time instances $U<T$, let us define event $A(U,T)$ as follows as the conjunction of the following:
\begin{itemize}
  \item[$(i)$] $X(U) \in [0,1)$,
  \item[$(ii)$] $X(t)$ is nonnegative at time $t$, for all $ U < t < T$, and
  \item[$(iii)$] $X(T) \ge \xi = \sqrt{n} \log n$.
\end{itemize}

Now we estimate the probability that $A(U,T)$ holds for some fixed $U,T$, $1 \le U \le T \le n$. Let $T_1$ be the smallest $t> U$ so that $X(t)<0$ or $t=T$. For all $t$ so that $U\leq t<T_1$,
$$
X(t+1)-X(t) =
\begin{cases}
1-c_1 &\mbox{\ with probability\ } (1-c_1-X(t)/t)^{1-\alpha} \le (1-c_1)^{1-\alpha} = c_1 \\
-c_1 &\mbox{\ otherwise.}
\end{cases}
$$
So for $U \le t \le T_1$
$$
\E(X(t+1)-X(t)~|~X(t)) ~~\le~~ (1-c_1)c_1 - c_1(1-c_1) ~~=~~ 0\,.
$$
Thus, $X(U),X(U+1),\dots,X(T_1)$ is a supermartingale where $|X(U)|\le 1$ and $|X(t+1)-X(t)| \le 1$. So from Lemma~\ref{lem:martingales} it
follows that
$$
\Prob \Big(X(T_1) \ge \xi \Big) \le \exp(-\xi^2/2n) ~~=~~ \exp(-\Theta(\log^2 n)) .
$$
Since condition ($ii$) is equivalent to $T_1=T$, we have $$\Prob(A(U,T))\leq\Prob(T_1=T)\Prob(X(T_1)\ge\xi |\,T=T_1),$$ and by the previous argument the last probability is exponentially small.

Similarly, we can define events $B(U,T)$, applying conditions $(i)$, $(ii)$ and $(iii)$ to $-X(t)$. (So condition $(i)$ of $B(U,T)$ is that $X(U)\in (-1,0]$, etc.) It can then be shown in an analogous way that $\Prob(B(U,T)) \le \exp(-\Theta(\log^2 n))$ for any $U,T,1 \le U \le T \le n$. Since all events have small probability \emph{wep} none of them occur. Indeed,
$$
\E \left( \sum_{U,T} I_{A(U,T)} + \sum_{U,T} I_{B(U,T)} \right) = O(n^2) \exp(-\Theta(\log^2 n)) = \exp(-\Theta(\log^2 n))
$$
and this fact follows from Markov's inequality. Given that none of the events occur, the assertion holds deterministically.
\end{proof}

We can repeat the same argument for all $Y_k$'s ($2 \le k \le K$). Since the error terms are cumulating, in order to get an asymptotic behaviour $K$ has to satisfy the following equation
$$
K \sqrt{n} \log n = K^{-1-1/\alpha} n / \log n = o(c_K n)\,.
$$
Note that $K(n)$ as defined earlier satisfies this equation. This completes the sketch of the proof of Theorem~\ref{thm:degree}.


\section{Deterministic Ranking Schemes}\label{sec:deterministic}

In this section we consider two ranking schemes that are deterministic, that is, the rank of a vertex $r(v_i,t)$ does not depend on $G_t$, but is completely determined by $i$ and $t$. In this case, the events that $v_i$ receives a link in time step $t$ are independent for all $t$. Thus, $\deg(v_i,n)$ is the sum of $n-i$ independent Bernouilli trials with pre-determined probabilities. The general theory about such sums can be directly applied to obtain the results in this section.

\bigskip

\subsection{Ranking by Age}\label{sec:age}

Ranking by age means that older vertices have a lower rank. Precisely, the rank of a vertex equals the time it is born, that is, $r(v_i,t) = i$ for all $1\leq i\leq t$. As mentioned in the Introduction, rank-based attachment with ranking by age is a special case of the growing protean graph model defined in~\cite{PPNW}. The growing protean graph model is more general since it permits deletion of vertices. Theorem~5.1 and Theorem~5.2 in~\cite{PPNW} give results for the degree of a vertex that apply to our model. However, the next theorem gives stronger concentration results since it is adopted to the special case that no deletion occurs.

\begin{theorem}\label{thm:age_degree}
For ranking by age, the expected degree of a vertex $v_i$, $i \in [n]$, is given by
$$
\E \deg(v_i,n) = \big(1+O(n^{-\alpha(1-\alpha)/3})\big) d \frac{1-\alpha}{\alpha} \left( \left(\frac ni\right)^\alpha+\frac {2\alpha-1}{1-\alpha} \right).
$$
Moreover, if $i<n/\log^{3/\alpha} n$, then \emph{wep}
$$
\deg(v_i,n) = \big(1+O(\log^{-1/2}n)\big) d \frac{1-\alpha}{\alpha} \left(\frac ni\right)^\alpha.
$$
\end{theorem}
\begin{proof}
Let $X(t,j)$ be a random indicator variable for an event that vertex $v_t$ joins $v_i$ at substep $j$ of step $t$ ($t\in [n]$, $j\in[d]$).
$$
\Prob(X(t,j)=1) = 1-\Prob(X(t,j)=0)=\
\begin{cases}
i^{-\alpha}/g_{\alpha}(t) & \textrm{\ for\ \ }  t > i\\
0 & \textrm{\ otherwise.}
\end{cases}
$$
The number of neighbours $v_t$ of $v_i$ such that $t>i$ is a random variable and can be expressed as a sum $\sum_{t=i+1}^n \sum_{j=1}^d X(t,j)$ of independent random variables. Since the number of neighbours $v_t$ of $v_i$ such that $t<i$ is always $d$,
$$
\deg(v_i,n) = d + \sum_{t=i+1}^n \sum_{j=1}^d X(t,j)\,.
$$
Thus, using~(\ref{denominator}),
$$
\E \deg(v_i,n) = d + d \sum_{t=i+1}^n \frac {i^{-\alpha}}{g_{\alpha}(t-1)} = d + d i^{-\alpha} \sum_{t=i+1}^n \frac {1-\alpha}{t^{1-\alpha}+O(1)}\,.
$$

Assuming that $i \ge n^{\alpha/3}$, we get that
\begin{eqnarray*}
  \E \deg(v_i,n) &=& d + \big(1+O(n^{\alpha(\alpha-1)/3})\big) d (1-\alpha) i^{-\alpha} \sum_{t=i+1}^n t^{\alpha-1} \\
  &=& d + \big(1+O(n^{\alpha(\alpha-1)/3})\big) d \frac{1-\alpha}{\alpha} i^{-\alpha} \left(n^\alpha-i^\alpha+O(i^{\alpha-1})\right) \\
  &=& \big(1+O(n^{\alpha(\alpha-1)/3})\big) d \frac{1-\alpha}{\alpha} \left( \left(\frac ni\right)^\alpha+\frac {2\alpha-1}{1-\alpha} \right)\,.
\end{eqnarray*}
A similar calculation can be done for $i < n^{\alpha/3}$, noting the fact that
$$
\E \deg(v_i,n) = O(n^{\alpha/3}) + d i^{-\alpha} \sum_{t=n^{\alpha/3}}^n \frac {1-\alpha}{t^{1-\alpha}+O(1)}\,.
$$

In order to finish the proof, we use the fact that a sum of independent random variables with large enough expected value is not too far from its mean (see, for example, Theorem~2.8 in~\cite{JLR}). From this it follows that, if $\eps \le 3/2$, then
\begin{equation}\label{eqn:Chernoff}
\Prob \left(|\deg(v_i,n)-\E \deg(v_i,n)| \ge \eps \E \deg(v_i,n)\right) \le 2 \exp \left(-\frac {\eps^2}{3} \E \deg(v_i,n) \right) .
\end{equation}
Note that $\E \deg(v_i,n) = \Omega(\log^3 n)$ for $i < n/\log^{3/\alpha} n$. If we let $\eps = \log n/\sqrt{\E \deg(v_i,n)}$ in~(\ref{eqn:Chernoff}), we get that \emph{wep}
$\deg(v_i,n) = \big(1+O(\eps)\big) \E \deg(v_i,n)$ and the assertion follows.
\end{proof}

Observe that, for small $i$, the expected degree of a vertex $v_i$ is dominated by the factor $d \frac{1-\alpha}{\alpha} \left(\frac ni\right)^\alpha$. Consequently, the degrees are distributed according to the power law. More specifically, let $Z_k=Z_k(n,d,\alpha)$ denote the number of vertices of degree $k$ and $Z_{\ge k}=\sum_{l\ge k} Z_l$. The following theorem shows that the $Z_{\ge k}$ follow a power law with exponent $1/\alpha$. Since the $Z_{\ge k}$ represent the cumulative degree distribution, this implies that the degree distribution follows a power law with exponent $1+1/\alpha$.

\begin{theorem}
Let $0<\alpha<1$ and $d \in \N$, $k \ge \log^4 n$. Then \emph{wep}
$$
Z_{\ge k} = \big(1+O(\log^{-1/3} n)\big) n \left(\frac {1-\alpha}{\alpha} \cdot \frac dk \right)^{1/\alpha}.
$$
\end{theorem}
\begin{proof}
This theorem is a simple consequence of Theorem~\ref{thm:age_degree}. One can show that \emph{wep} each vertex $v_i$ such that
$$
i \ge \big(1+\log^{-1/3} n \big) n \left(\frac {1-\alpha}{\alpha} \cdot \frac dk \right)^{1/\alpha}
$$
has fewer than $k$ neighbours, and each vertex $v_i$ for which
$$
i \le \big(1-\log^{-1/3} n \big) n \left(\frac {1-\alpha}{\alpha} \cdot \frac dk \right)^{1/\alpha}
$$
has more than $k$ neighbours.
\end{proof}

\bigskip

\subsection{Ranking by inverse age}\label{sec:inverse-age}

To contrast the other schemes, we considered a scheme where new vertices are ranked the highest. Precisely, $r(v_i,t)=t-i+1$. Intuitively, this scheme breaks the effect of ``cumulative advantage'', since no vertex has high rank long enough to accumulate a high degree. The results from this section give evidence that, indeed, this scheme does  not lead to a power law degree distribution.

Note that
\begin{eqnarray*}
\E \deg(v_1,n) &=& d + d \sum_{t=2}^n \frac {(t-1)^{-\alpha}} {g_\alpha(t-1)} ~ ~ = ~ ~(1+o(1)) d (1-\alpha) \sum_{t=1}^{n-1} \frac {1}{t}\\
&=& (1+o(1)) d (1-\alpha) \log n
\end{eqnarray*}
and also it is not hard to see that $\E \deg(v_i,n) > \E \deg(v_j,n)$ for $1 \le i < j \le n$. Thus $\E \deg(v_i,n) < (1+o(1)) d (1-\alpha) \log n$ for all $i \in [n]$. Again $\deg(v_i,n)$ can be expressed as a sum of independent $0-1$ random variables but since the expected degree is so low we cannot hope for concentration; the Chernoff bound only tells us that \emph{wep} the maximum degree of $G_n$ is $\E \deg(v_1,n) + O(\log {n})$.

We also show that the number of vertices with expected degree at least $k$ decreases exponentially with $k$. This suggests that the degree distribution does not follow a power law.

\begin{theorem} \label{thm:inverse_age}
Let $0<\alpha<1$, $d \in \N$, and $i = i(n) \in [n]$. The expected degree of a vertex $v_i$ satisfies the following inequalities
\begin{eqnarray*}
\E \deg(v_i,n) &\ge& d + (1+o(1)) d (1-\alpha) \alpha \log (n-i)\\
\E \deg(v_i,n) &\le& d + (1+o(1)) d (1-\alpha) \log (n-i) \,.
\end{eqnarray*}
\end{theorem}
\begin{proof}
Define
$$
f(i) = \sum_{t = 0}^{n-i-1} \frac {1}{(t+1)^\alpha (t+i)^{1-\alpha}}.
$$
Then, using the approach as in Theorem~\ref{thm:age_degree}, we obtain,
\begin{eqnarray*}
\E \deg(v_i,n) &=& d + d \sum_{t=i+1}^n \frac {(t-i)^{-\alpha}}{g_\alpha(t-1)}\\
&=& d + (1+o(1)) d (1-\alpha) \sum_{t=i+1}^n \frac {(t-i)^{-\alpha}}{(t-1)^{1-\alpha}}\\
&=& d + (1+o(1)) d (1-\alpha) f(i)\,,
\end{eqnarray*}
for any $i \in [n]$. The assertion follows from the fact that
\begin{eqnarray*}
f(i) &\ge& \int_{0}^{n-i} \frac {1}{(x+1)^\alpha (x+i)^{1-\alpha}} dx
~~=~~ \int_{0}^{n-i} \left(1+\frac {i-1}{x+1} \right)^\alpha \frac {1} {x+i} dx \\
&\ge& \int_{0}^{n-i} \left(1+\alpha \frac {i-1}{x+1} \right) \frac {1}
{x+i} dx ~~=~~ \int_{0}^{n-i} \left(\frac {1-\alpha}{x+i} + \frac {\alpha}{x+1} \right) dx \\
&=& (1-\alpha) \log (n/i) + \alpha \log (n-i+1) ~~\ge~~ \alpha \log (n-i)
\end{eqnarray*}
and
\begin{eqnarray*}
f(i) &\le& 1 + \int_{0}^{n-i-1} \frac {1}{(x+1)^\alpha (x+i)^{1-\alpha}}
dx ~~=~~ 1 + \int_{0}^{n-i-1} \left(\frac {x+1}{x+i}
\right)^{1-\alpha} \frac {1}{x+1} dx \\
&\le& 1 + \left(\frac {n-i}{n-1} \right)^{1-\alpha} \log(n-i) ~~\le~~
\log (n-i) +1 \,.
\end{eqnarray*}

\begin{cor}
Let $0<\alpha<1$, $d \in \N$.
\end{cor}
\begin{eqnarray*}
\# \{v_i : \E \deg(v_i,n) \ge k\} &\ge& n - (1+o(1)) \exp \left(\frac {k-d}{d(1-\alpha)\alpha} \right) \\
\# \{v_i : \E \deg(v_i,n) \ge k\} &\le& n - (1+o(1)) \exp
\left(\frac {k-d}{d(1-\alpha)} \right) \,.
\end{eqnarray*}
\end{proof}


\section{Random ranking}

In the two ranking schemes discussed in this section, the initial rank $r(v_t,t)$ of a new vertex $v_t$ is a random variable $R_t\in [t]$. The new rank function is simply formed by inserting the new vertex into the existing ranking, so for all $j\in[t-1]$, $r(v_j,t)=r(v_j,t-1)$ if $r(v_j,t-1)<R_t$, and $r(v_j,t)=r(v_j,t-1)+1$ otherwise. The difference in the two schemes lies in the way that $R_t$ is chosen: in the first scheme, the rank of each vertex is based on a fixed, but randomly chosen, label, while in the second scheme, $R_t$ is randomly drawn from $[t]$.

\subsection{Ranking by random labeling}\label{sec:labels}

In this scheme, each new vertex $v_t$ obtains a label $l(v_t)\in (0,1)$ chosen uniformly at random. (Note that the probability that two vertices receive the same label is zero.) Vertices are ranked by their labels: if $l(v_i)<l(v_j)$, then $r(v_i,t)<r(v_j,t)$.

First we note that the process of choosing a label \emph{uar} from $(0,1)$ does not imply loss of generality. Namely, suppose that the labels are chosen from $\R$ according to any probability distribution with a strictly increasing {\sl cumulative distribution function} $F$. Since $F$ is an increasing function, labels $F(l(v_i))$ lead to exactly the same ranking as labels $l(v_i)$. But $\pp(F(l(v_i))\leq x)=\pp(l(v_i)\leq F^{-1}(x))=F(F^{-1}(x))=x$, so the values of labels $F(l(v_i))$ are chosen from $(0,1)$ according to the uniform distribution.

Since the arguments used here are largely similar to those used in Section~\ref{sec:deterministic} (namely, applying the Chernoff bounds for a sum of independent random variables), we omit technical details.

\begin{theorem} \label{thm:labels_degree}
Let $0<\alpha<1$, $d \in \N$, $i = i(n) \in [n]$, and $0 < l(v_i)=l(v_i)(n) < 1$. If $i \cdot l(v_i) > \log^3 n$, then the expected degree of a vertex $v_i$ that obtained a label $l(v_i)$, is given by
$$
\E \deg(v_i,n) = d + (1+O(\log^{-(1-\alpha)/2} n)) d (1-\alpha) l(v_i)^{-\alpha} \log(n/i) \,,
$$
and \emph{wep}
$$
\deg(v_i,n)= \E \deg(v_i,n) + O(\sqrt{\E \deg(v_i,n)} \log n) \,.
$$
\end{theorem}
\begin{proof}
Note that $r(v_i,t)$ is the sum of independent indicator variables of the events $l(v_j)\leq l(v_i)$ for $i<j\leq t$. Using large deviation inequalities and the fact that $i \cdot l(v_i) > \log^3 n$, we get that \emph{wep} for all $i<t\leq n$,
$$
t \cdot l(v_i) (1-\log^{-1/2} n) \le r(v_i,t) \le t \cdot l(v_i) (1+\log^{-1/2} n)\,.
$$
Thus,
\begin{eqnarray*}
\E \deg(v_i,n) &=& d + d \sum_{t=i+1}^n \frac {\big(t \cdot l(v_i) (1+O(\log^{-1/2} n))\big)^{-\alpha}}{g_\alpha(t-1)} \\
&=& d + (1+O(\log^{-(1-\alpha)/2} n)) d (1-\alpha) l(v_i)^{-\alpha} \sum_{t=i+1}^n \frac {1}{t} \\
&=& d + (1+O(\log^{-(1-\alpha)/2} n)) d (1-\alpha) l(v_i)^{-\alpha} \log(n/i)\,.
\end{eqnarray*}
Since $\deg(v_i,n)$ can be expressed as a sum of independent random variables, we can use the Chernoff bound to show the concentration result.
\end{proof}

Using the notation from Section~\ref{sec:deterministic} we present the main result.

\begin{theorem}\label{thm:labels_powerlaw}
Let $0<\alpha<1$ and $d \in \N$, $\log^3 n \le k \le n^{\alpha/2}/\log^{3\alpha} n$. Then \emph{wep}
$$
Z_{\ge k} = (1+o(1)) n \left( \frac{d(1-\alpha)}{k} \right)^{1/\alpha} \Gamma \left( \frac {1}{\alpha} +1 \right).
$$
\end{theorem}
\begin{proof}
From Theorem~\ref{thm:labels_degree} and the fact that $k \ge \log^3 n$ it follows that \emph{wep} all vertices $v_i$ such that $i \ge i_k = k^{1/\alpha} \log^4 n$ and
$$
l(v_i) \ge (1+\log^{-(1-\alpha)/3} n) \left( \frac{d (1-\alpha) \log(n/i)}{k} \right)^{1/\alpha}
$$
has fewer than $k-d=(1+o(1))k$ neighbours, and each vertex $v_i$ for which
$$
l(v_i) \le (1-\log^{-(1-\alpha)/3} n) \left( \frac{d (1-\alpha) \log(n/i)}{k} \right)^{1/\alpha}
$$
has more than $k$ neighbours.

Let $X_i^+$, $X_i^-$, $i \in [n]$, be a family of independent random variables such that
$$
\Prob(X_i^+=1)=1-\Prob(X_i^+=0)=
\begin{cases}
(1+\log^{-(1-\alpha)/3} n) \left(  \frac{d (1-\alpha) \log(n/i)}{k} \right)^{1/\alpha} & \textrm{\ for\ \ }  i \ge i_k\\
1 & \textrm{\ otherwise,}
\end{cases}
$$
and
$$
\Prob(X_i^-=1)=1-\Prob(X_i^-=0)=
\begin{cases}
(1-\log^{-(1-\alpha)/3} n) \left(  \frac{d (1-\alpha) \log(n/i)}{k} \right)^{1/\alpha} & \textrm{\ for\ \ }  i \ge i_k\\
0 & \textrm{\ otherwise.}
\end{cases}
$$
Then, $Z_{\ge k}$ can be bounded from below by $X^-=\sum_{i=1}^n X_i^-$ and from above by $X^+=\sum_{i=1}^n X_i^+$. Thus,
\begin{eqnarray*}
\E Z_{\ge k} &=& O(k^{1/\alpha} \log^4 n) + (1+o(1)) \sum_{i=1}^n \left( \frac{d (1-\alpha) \log(n/i)}{k} \right)^{1/\alpha}\\
&=& O(k^{1/\alpha} \log^4 n) + (1+o(1)) \left( \frac{d(1-\alpha)}{k} \right)^{1/\alpha} n \int_0^1 \left( \log(1/x) \right)^{1/\alpha} dx \,,
\end{eqnarray*}
and putting $u = \log(1/x)$ we get
\begin{eqnarray*}
\E Z_{\ge k} &=& O(k^{1/\alpha} \log^4 n) + (1+o(1)) \left( \frac{d(1-\alpha)}{k} \right)^{1/\alpha} n \int_0^{\infty} u^{1/\alpha} e^{-u} du \\
&=& (1+o(1)) n \left( \frac{d(1-\alpha)}{k} \right)^{1/\alpha} \Gamma \left( \frac {1}{\alpha} +1 \right) \,,
\end{eqnarray*}
where $\Gamma(x)$ denotes the (complete) gamma function. Since the gamma function is an extension of the factorial, and is increasing for $x\geq 2$, $\Gamma(1/\alpha+1)$ is a constant which lies between $\lfloor 1/\alpha \rfloor!$ and $\lceil 1/\alpha \rceil!$.

Since $k \le n^{\alpha/2}/\log^{3\alpha} n$, $\E Z_{\ge k} = \Omega(\sqrt{n} \log^3 n)$. Using large deviation inequalities one more time, we can show that \emph{whp} $Z_{\ge k} = (1+o(1)) \E Z_{\ge k}$. This finishes the proof of the theorem.
\end{proof}

\bigskip

\subsection{Randomly chosen initial rank}\label{sec:random}

Next, we consider the case where the rank of the new vertex $v_t$, $R_t=r(v_t,t)$, is chosen at random from $[t]$. As described earlier, the ranks of existing vertices are adjusted accordingly. In contrast to the previous scheme, in this case it now matters according to which distribution $R_t$ is chosen. We make the assumption that all initial ranks are chosen according to a similar distribution. In particular, we fix a continuous bijective function $F:[0,1]\rightarrow[0,1]$, and for all integers $1\leq k\leq t$, we let
\begin{equation*}\label{F}
\Prob(R_t\leq k)=F\left(\frac kt\right)\,.
\end{equation*}

Thus, $F$ represents the limit, for $t$ going to infinity, of the cumulative distribution functions of the variables $R_t$. To simplify the calculations while exploring a wide array of possibilities for $F$, we assume $F$ to be of the form
\[
F(x)=x^s\mbox{, where }s>0.
\]
A special case is the case $s=1$, where the distribution of each $R_t$ is uniform. We will see that this case is the threshold for a power law degree distribution; if $s<1$, then the probability that a new vertex receives a high rank (that is, a low value of $R_t$) is enhanced, and thus we get behaviour similar to that of  age-based ranking, as seen in Subsection~\ref{sec:age}, including a power law degree distribution; if $s>1$, then the probability that a new vertex receives a high rank is diminished, and we get behaviour similar to the inverse age ranking scheme described in Subsection~\ref{sec:inverse-age}.

To study the degree of a given vertex $v_i$ in $G_n$ under this ranking scheme we again use the differential equations method. Assume that vertex $v_i$ obtained an initial rank $R_i$. Then $r(v_i,t)$, $t > i$, is a random variable, which in time step $t$ increases by one precisely when $R_{t} \le r(v_i,t-1)$. Since the latter happens with probability $F(r(v_i,t-1)/t)$, we have that
\begin{equation}
\label{cond_exp_random_rank}
\E (r(v_i,t)-r(v_i,t-1) ~|~ G_t) = \left(\frac {r(v_i,t-1)}{t}\right)^s\,.
\end{equation}
Using a real function $z(x)$ to model the behaviour of $r(v_i,xn)/n$, the above equation leads to the following differential equation for $z$:
\begin{equation}\label{z:de}
z'(x) = \left(\frac {z(x)}{x}\right)^s
\end{equation}
with the initial condition $z(i/n)=R_i/n$.

If $s=1$, the general solution is $z(x)=C x$, $C \in \R$ and the particular solution is $z(x) = \frac{R_i}i x$. This \emph{suggests} that a random variable $r(v_i, t)$ should be close to a deterministic function $\frac{R_i}i t$. We will use martingales to show that this is indeed the case.

Let $Y_t=\frac{r(v_i,t)}{t+1}$ for all $i\leq t\leq n$. The sequence $\{Y_t:i\leq t\leq n\}$ is a martingale with respect to the random process $\{G_t\}$. Namely,
\begin{eqnarray*}
\E (Y_t ~|~ Y_{t-1})
&=& \frac {r(v_i,t-1)+\left(\frac {r(v_i,t-1)}{t}\right)}{t+1} ~=~ \frac {r(v_i,t-1)}{t} ~=~ Y_{t-1}\,.
\end{eqnarray*}
In order to show a concentration for $Y_t$, and thus for $r(v_i,t)$, we use a well-known Azuma-Hoeffding inequality (see for example Lemma~4.1 in~\cite{NW}).

\begin{lemma}\label{lem:martingales2}
Let $X_0, X_1, \dots, X_t$ be a martingale such that $|X_j-X_{j-1}| \le c_j$, $1 \le j \le t$, for constants $c_j$, Then for any $\alpha>0$
$$
\Prob(|X_t-X_0|\ge \alpha) \le 2 \exp \left(- \ \frac {\alpha^2}{2 \sum c_j^2} \right)\,.
$$
\end{lemma}

Now, we are ready to state a concentration theorem.
\begin{theorem}
Let $i=i(n) \in [n]$, $R_i=R_i(n)$ such that $1 \le  R_i\le i$ and $R_i^2/(i+1) > \log^4 n$. If the vertex $v_i$ obtained an initial rank $R_i$, and $l(v_i)=R_i/(i+1)$, then \emph{wep}
$$
t \cdot l(v_i) (1-\log^{-1/2} n) \le r(v_i,t) \le t \cdot l(v_i) (1+\log^{-1/2} n)
$$
for all $t, i<t\le n$.
\end{theorem}
\begin{proof}
Let $Y_t$ be a random variable defined as before. Note that $l(v_i)=Y_i$. It has been shown that $\{Y_t\}$ is a martingale and it is also easy to see that $|Y_t-Y_{t-1}| \le \frac {1}{t}$. Since
$$
\sum_{t=i}^n \left(\frac 1t \right)^2 = \int_{i}^n x^{-2} dx + O(1) = O \left(\frac 1i \right)\,,
$$
we can apply Lemma~\ref{lem:martingales} with $c_t=\frac 1t$ and $\alpha = \frac {\log n}{\sqrt{i}}$ to obtain that \emph{wep} $|Y_t - Y_i| \le \alpha$. So \emph{wep}
$$
Y_t = Y_i \left(1 + O \left(\frac {\alpha}{Y_i} \right) \right) = l(v_i) \left(1 + O\left(\frac {\log n}{\sqrt {i} \cdot l(v_i)}\right) \right) = l(v_i) (1 + O(\log^{-1} n))\,.
$$
Moreover, one can use a concept of a stopping time (see Section~12.4 in~\cite{GS}) to show that \emph{wep} all values of $Y_t$, $i\leq t\leq n$, lie within the bounds given by the equation above. This finishes the proof.
\end{proof}

Since the proofs of Theorems~\ref{thm:random_degree} and~\ref{thm:random_powerlaw} are almost the same as the proofs of Theorems~\ref{thm:labels_degree} and~\ref{thm:labels_powerlaw}, we omit them stating the results only.

\begin{theorem} \label{thm:random_degree}
Let $0<\alpha<1$, $d \in \N$, $i = i(n) \in [n]$, $R_i=R_i(n)$ such that
$1 \le  R_i\le i$. If the vertex $v_i$ obtained an initial rank $R_i$ such that $R_i^2/(i+1) > \log^4 n$, and $l(v_i)=R_i/(i+1)$, then the expected degree of  $v_i$ is
$$
\E \deg(v_i,n) = d + (1+O(\log^{-(1-\alpha)/2} n)) d (1-\alpha) l(v_i)^{-\alpha} \log(n/i) \,.
$$
and \emph{wep}
$$
\deg(v_i,n)= \E \deg(v_i,n) + O(\sqrt{\E \deg(v_i,n)} \log n) \,.
$$
\end{theorem}

\begin{theorem}\label{thm:random_powerlaw}
Let $0<\alpha<1$ and $d \in \N$, $\log^3 n \le k \le n^{\alpha/3}/\log^{2\alpha} n$. Then \emph{wep}
$$
Z_{\ge k} = (1+o(1)) n \left( \frac{d(1-\alpha)}{k} \right)^{1/\alpha} \Gamma \left( \frac {1}{\alpha} +1 \right).
$$
\end{theorem}

Note that Theorem~\ref{thm:labels_powerlaw} and~\ref{thm:random_powerlaw} suggest that the random ranking scheme with uniform distribution behaves similarly to the random labeling scheme. However, the upper bound values $k$ for which the power law holds is lower in this case. This can be explained by the fact that the eventual rank of a vertex is not always easy to predict in this scheme. For example, assume that in the random labeling scheme, the first vertex obtained a label of $1/2$. Then its rank at time $n$ is almost surely close to its expected value of $n/2$. In the current scheme, if a vertex has initial rank $R_i/2i$, then the expected rank at time $n$ also equals $n/2$ but this rank is not concentrated. Namely, the rank behaves like the proportion of white balls in Polya's urn problem, and thus $r(v_1,n)/n$ converges to a random variable with uniform distribution on $[0,1]$.

\bigskip

Next, we consider the case where $s>1$. The general solution of the differential equation~(\ref{z:de}) is $z(x)=(x^{1-s}+C)^{-1/(s-1)}$. Thus, we expect $r(v_i,t)$ to be approximately equal to $nz(t/n)=(t^{1-s}+c)^{-1/s-1}$. Note that, if $t$ gets large, this function converges to a constant. As we will see, the definition below captures the value of this constant.
\begin{equation} \label{defRi}
R^*_i= \left(R_i^{1-s}-(i+1)^{1-s}\right)^{-1/(s-1)}.
\end{equation}

\begin{theorem}
\label{thm:rank}
For all $i\geq n^{1/2}\log^{s+1}n$, if vertex $v_i$ has initial rank $R_i$ so that $n^{1/2}\log^{s+1}{n}\leq R_i\leq (1-\log^{1-s}{n})i$, then \emph{wep}
\[
r(v_i,t)=R^*_i \left(1+\left(\frac{R_i^*}{t+1}\right)^{s-1}\right)^{\frac{-1}{s-1}}(1+O(\log^{-1} n)).
\]
\end{theorem}

The proof uses the supermartingale method of Pittel et al.~\cite{PSW}, as described in~\cite[Corollary 4.1]{NW}. We  need the following lemma.

\begin{lemma} \label{l:corollary}
Let $G_{0},G_{1},\dots ,G_{n}$ be a random process and $ X_{t}$ a random variable determined by $G_{0},G_{1},\dots ,G_{t}$, $0\leq t\leq n$. Suppose that for some real $\beta$ and constants $\gamma_{t}$,
\begin{equation*}
\E (X_{t}-X_{t-1}~|~G_{0},G_{1},\dots ,G_{t-1})<\beta
\end{equation*}
and
\begin{equation*}
|X_{t}-X_{t-1}-\beta|\leq \gamma_{t}
\end{equation*}
for $1\leq t\leq n$. Then for all $\alpha >0$,
\[
\Prob \big(\mbox{For some $t$ with }0\leq t\leq n:X_{t}-X_{0}\geq t\beta+\alpha \big)\leq \exp \Big(-\frac{\alpha ^{2}}{2\sum_{j=1}^n \gamma_{t}^{2}}\Big)\;.
\]
\end{lemma}

\noindent
{\sl Proof of Theorem \ref{thm:rank}}.
We transform $r(v_i,t)$ into something close to a martingale. Consider the following real-valued function
\begin{equation}
H(x,y)=y^{1-s}-(x+1)^{1-s}  \label{e:H}
\end{equation}
Let $\mathbf{w}_{t}=(t, r(v_i,t))$, and consider the sequence of random variables $(H(\mathbf{ w}_{t}):i\leq t\leq n)$. Note that $H(i,R_i)=(R^*_i)^{1-s}$. We will show that \emph{wep} $H(t,r(v_i,t))$ is close to $H(i,R_i)$. The function $H$ is chosen so that $H(\mathbf{w})$ is constant along every trajectory $\mathbf{w}$ of the differential equation~(\ref{z:de}).

Note that
$$
\mathrm{grad}\mbox{ }H(\mathbf{w}_{t})=\left(-(1-s)(t+1)^{-s},(1-s)r(v_i,t)^{-s}\right) ,
$$
It follows from the choice of $H$, and can be checked using~(\ref{cond_exp_random_rank}), that
$$
\mathbb{E}(\mathbf{w}_{t+1}-\mathbf{w}_{t}~|~G_{t})\cdot \mbox{\ \textrm{grad} }H(\mathbf{w}_{t})=0,
$$

Using the fact that $R_i\leq r(v_i,t)\leq t$ at all times, we can show that all second-order partial derivatives of $H$ evaluated at $\mathbf{w}_t$ are $O(R_i^{-(s+1)})$. Therefore,
\begin{equation}\label{e:h_diff}
H(\mathbf{w}_{t+1})-H(\mathbf{w}_{t})=(\mathbf{w}_{t+1}-\mathbf{w}_{t})\cdot \mathrm{grad}\mbox{ }H(\mathbf{w}_{t})+O(R_i^{-(s+1)}).
\end{equation}

Taking the expectation of~(\ref{e:h_diff}) conditional on $G_{t}$, we obtain that
$$
{\mathbb{E}}(H(\mathbf{w}_{t+1})-H(\mathbf{w}_{t})~|~G_{t})=O(R_i^{-(s+1)}).
$$

The rank changes by at most one in each step, so from the above, we obtain
\begin{eqnarray*}
|H(\mathbf{w}_{t+1})-H(\mathbf{w}_{t})| &\leq& (s-1) \left(r(v_i,t)^{-s}+(t+1)^{-s}\right)+O(R_i^{-(s+1)}) \\
&=&O(R_i^{-s}).
\end{eqnarray*}

Now we may apply Lemma~\ref{l:corollary} to the sequence $(H(\mathbf{w} _{t}):i\leq t\leq n)$, and symmetrically to $(-H(\mathbf{w}_{t}):i\leq t\leq n)$, with $\alpha =R_i^{-s}n^{1/2}\log{n}$, $\beta=O(R_i^{-(s+1)})$, and $\gamma_{t}=O(R_i^{-s})$. From the lower bound on $R_i$ it follows that $n\beta=O(\alpha)$, and we obtain that \emph{wep}
$$
|H(\mathbf{w}_{t})-H(\mathbf{w}_{i})|=O(R_i^{-s} n^{1/2}\log{n})
$$
for $i \le t \le n$. As $H(\mathbf{w}_{i})=(R_i^*)^{1-s}$, this implies from the definition~(\ref{e:H}) of the function $H$, that \emph{wep}
\begin{eqnarray*}
r(v_i,t)^{1-s}&=&(R^*_i)^{1-s}+(t+1)^{1-s}+O(R_i^{-s}n^{1/2}\log{n})\\
&=& \left( (R^*_i)^{1-s}+(t+1)^{1-s}\right)(1+O((R^*_i)^{s-1} R_i^{-s}n^{1/2}\log{n})
\end{eqnarray*}
for $i \le t \le n$, so
\begin{equation*}
r(v_i,t)=R_i^* \left(1+\left(\frac{R_i^*}{t+1}\right)^{s-1}\right)^{\frac{-1}{s-1}}(1+O((R^*_i)^{s-1} R_i^{-s} n^{1/2}\log{n} )).
\end{equation*}
Since $R_i/i\leq (1-\log^{1-s}{n})$, we have that
\begin{equation} \label{ratio}
R_i^*/R_i=O(\log{n}).
\end{equation}
Since  $R_i\geq n^{1/2} \log^{s+1}{n}$ we have that
$$
(R^*_i)^{s-1} R_i^{-s} n^{1/2}\log{n} = O\left( R_i^{-1} n^{1/2} \log^{s} {n} \right) = O(\log^{-1} n)
$$
which finishes the proof of the theorem.\hfill $\Box$

We can now use the same approach as for age-based ranking.

\begin{theorem}
\label{thm:deg}
For a vertex $v_i$ so that $n^{1/2}\log^{s+1}{n}\leq R_i\leq (1-\log^{1-s}{n})i$ and $R_i^*\leq n\log^{-3/\alpha}n$, \emph{wep},
\[
\deg(v_i,n)= (1+O(\log^{-\min\{1/2,3(s-1)/\alpha\}}{n}))d\frac{a}{1-\alpha}\left(\frac{n}{R_i^*}\right)^\alpha.
\]
Moreover, for all vertices $v_i$, \emph{wep}
\[
\deg(v_i,n)= d \frac{1-\alpha}{\alpha} \left(\frac{n}{R_i}\right)^\alpha + O(\log^2 n).
\]
\end{theorem}

\begin{proof}
The proof follows the same reasoning as the proof of Theorem~\ref{thm:age_degree}. To prove the first part, using Theorem~\ref{thm:rank}, we obtain the expected degree of $v_i$ at time $n$ as follows:
\begin{eqnarray*}
\E \deg(v_i,n) &=& \E \deg(v_i,2R_i^*)+ d \sum_{t=2R_i^*+1}^n \frac{r(v_i,t)^{-\alpha}}{g_\alpha(t)}.
\end{eqnarray*}

For the first term, we use the fact that $r(v_i,t)\geq R_i$ for all $t\geq i$.
\begin{eqnarray}
\label{eq:Ri}
\E \deg(v_i,2R_i^*)&\leq &d+d\sum_{t=i+1}^{2R_i^*} \frac{R_i^{-\alpha}}{g_\alpha(t)}\\
&= & d + (1+o(1)) d R_i^{-\alpha} (1-\alpha) \sum_{t=i+1}^{2R_i^*} t^{\alpha -1}\notag\\
&=& O((R_i^*/R_i)^\alpha).\notag
\end{eqnarray}
Since $(R_i^*/R_i)= O(\log n)$ (see~(\ref{ratio})) and, by assumption, $n/R_i^*\geq \log^{3/\alpha} n$, we have that
\[
\E \deg(v_i,2R_i^*)= O(\log^\alpha n) = O(\left(n/R_i^*\right)^\alpha\log^{-2}n).
\]

Now, we can estimate the second part as follows:
\begin{align*}
d & \sum_{t=2R_i^*+1}^n \frac{r(v_i,t)^{-\alpha}}{g_\alpha(t)}  \\
&= (1+O(\log^{-1}n)) d (1-\alpha) (R_i^*)^{-\alpha} \sum_{t=2R_i^*+1}^n \frac{\left(1+\left(\frac{R_i^*}{t+1}\right)^{s-1}\right)^{\frac{\alpha}{s-1}}}{t^{1-\alpha}} \\
&= (1+O(\log^{-1}n)) d (1-\alpha) (R_i^*)^{-\alpha} \sum_{t=2R_i^*+1}^n \left( t^{\alpha-1} + O((R_i^*)^{s-1} t^{\alpha-s}) \right) \\
&= (1+O(\log^{-1}n)) d (1-\alpha) (R_i^*)^{-\alpha} \left( n^{\alpha}/\alpha - O((2R_i^*)^{\alpha}) + O((R_i^*)^{s-1} n^{\alpha-s+1}) \right) \\
&= (1+O(\log^{-\min\{1,3(s-1)/\alpha\}} n)) d \frac{1-\alpha}{\alpha} \left(\frac {n}{R_i^*}\right)^{\alpha}
\end{align*}
since $R_i^*\leq n \log^{-3/\alpha}n$. Therefore,
$$
\E  \deg(v_i,n) = (1+O(\log^{-\min\{1,3(s-1)/\alpha\}} n)) d \frac{1-\alpha}{\alpha} \left(\frac {n}{R_i^*}\right)^{\alpha}.
$$

Using the Chernoff bound as before (see~(\ref{eqn:Chernoff})), together with the fact that $\E \deg(v_i) = \Omega(\log^3 n)$ for $R_i^* < n\log^{-3/\alpha} n$, and putting $\eps = \log n/\sqrt{\E \deg(v_i)}$ in~(\ref{eqn:Chernoff}), we get that \emph{wep} $\deg(v_i) = \big(1+O(\eps)\big) \E \deg(v_i)$ and the assertion follows.

To prove the second part, we can use a calculation similar to~(\ref{eq:Ri}) to show that $\E \deg(v_i,n) = d\left(\frac{1-\alpha}{\alpha}\right)\left(\frac{n}{R_i}\right)^\alpha+O(1)$, and use the Chernoff bound to prove the statement of the theorem.
\end{proof}

\begin{theorem}
Let $k$ be so that $\log^{4}n\leq k\leq n^{\alpha/2}\log^{-\alpha(s+3)}n$. For random ranking with initial rank distribution given by $F$ where $F(x)=x^s$ and $s>1$,
$$
Z_{\ge k} = \big(1+O(\log^{1-s} n)\big) n \left(\frac {1-\alpha}{\alpha}
\cdot \frac dk \right)^{1/\alpha},
$$
\end{theorem}
\begin{proof}
Let $\omega(n)=n^{1/2}\log^{s+1}n$. Fix $k$ so that $\log^{\max(3\alpha,3/\alpha +2)} \leq k\leq n^{\alpha/2}\log^{-\alpha(s+3)}n$. Define sets $S_k^+$ and $S_k^-$ as follows:
\begin{eqnarray*}
S_k^-&=&\left\{ v_i\,|\,R_i^* \le \big(1-\log^{1-s} n \big) n \left(\frac {1-\alpha}{\alpha} \cdot \frac dk \right)^{1/\alpha}\right\},\\
S_k^+&=&\left\{ v_i\,|\,R_i^* \ge \big(1+\log^{1-s} n \big) n \left(\frac {1-\alpha}{\alpha} \cdot \frac dk \right)^{1/\alpha}\right\}.
\end{eqnarray*}
We will argue below that \emph{wep} all but a small fraction of the vertices in $S_k^-$ have degree at least $k$, and in $S_k^+$ have degree less than $k$. First, we estimate the size of $S_k^-$ and $S_k^+$.

Let $f(k)$ be a function so that $f(k)=\Theta(k^{-1/\alpha}n)$; $f(k)$ is meant to represent the bound on $R_i^*$ that defines $S_k^+$
or $S_k^-$. The bounds on $k$ imply that $f(k)=O(n\log^{-3}n)$ and $f(k)=\Omega(\omega(n)\log^2(n))$.

From (\ref{defRi}), $R_i^*\leq f(k)$ if and only if 
$R_i\leq ((i+1)^{1-s}+f(k)^{1-s})^{-1/(s-1)}$. 
Thus for any $i$, the probability that $R_i^*\leq f(k)$ equals $\left(1+\left(\frac{f(k)}{i+1}\right)^{1-s}\right)^{\frac{-s}{s-1}}(1+\frac1i)$. The expected number of vertices $v_i$ so that $R_i^*\leq f(k)$ is expressed by the following sum:
\begin{eqnarray*}
\sum_{i=1}^n \left(1+\left(\frac{f(k)}{i+1}\right)^{1-s}\right)^{\frac{-s}{s-1}}(1+\frac1i)
&=& \int_1^n \left(1+\left(\frac{f(k)}{x}\right)^{1-s}\right)^{\frac{-s}{s-1}}dx +O(\log n)\\
&=& f(k)\int_{1/f(k)}^{n/f(k)} (1+y^{s-1})^{\frac{-s}{s-1}}dy +O(\log n)\\
&=& f(k)(1 +O(n^{-1/2})).
\end{eqnarray*}
The last step can be explained as follows. The antiderivative of $(1+y^{s-1})^{\frac{-s}{s-1}}$ equals $(1+y^{1-s})^{\frac{1}{1-s}}$, and thus $\int_0^{\infty} (1+y^{s-1})^{\frac{-s}{s-1}}dy=1$. The integral from 0 to $1/f(k)$ is at most $1/f(k)=O(n^{-1/2})$. The integral from $n/f(k)$
to infinity equals $\frac{1}{s-1}\left(\frac{f(k)}{n}\right)^{s-1}(1+o(1))=o(1)$. The result then follows because $f(k)=\Omega(n^{1/2}\log n)$.

Using the Chernoff bound, and the lower bound on $f(k)$, it follows that \emph{wep}, the number of vertices with $R_i^*\leq f(k)$ equals
$(1+O(\log^{-1}n) f(k)$. Therefore, \emph{wep}
\[
|S_k^-|=\big(1-\log^{1-s} n \big) n \left(\frac {1-\alpha}{\alpha} \cdot \frac dk \right)^{1/\alpha},
\]
while the number of vertices that is neither in $S_k^+$ nor in $S_k^-$ is $O(\log^{1-s}|S_k^-|)$.

Consider the vertices in $S_k^-$. Let $f(k)= \big(1-\log^{1-s} n \big) n \left(\frac {1-\alpha}{\alpha} \cdot \frac dk \right)^{1/\alpha}$. From the bounds on $k$ it follows that $R_i^*\leq f(k)=O(n\log^{-3}n)$, so we may assume that $R_i^*\leq n\log^{-3/\alpha}n$. By Theorem~\ref{thm:deg}, if $\omega(n)\leq R_i\leq (1-\log^{1-s}n)i$, then \emph{wep} $\deg(v_i,n)\geq k$.

Consider the vertices in $S_k^-$ with initial rank $R_i<\omega(n)$. Since lower initial rank will \emph{wep} lead to higher degree, and since $\deg(v_i,n)$ would have been at least $k$ even if the initial rank $R_i$ had been $2\omega(n)$, we can conclude that these vertices also have degree at least $k$.

If $v_i\in S_k^-$ and $R_i>(1-\log^{1-s}n)i$, then this implies that $i\leq i_k$, where $i_k$ is so that
\begin{eqnarray*}
(1-\log^{1-s}n)i_k &=& ((i_k+1)^{1-s}+f(k)^{1-s})^{-1/(s-1)}\\
&=&\left(1+\left(\frac{f(k)}{i_k+1}\right)^{1-s}\right)^{-1/(s-1)}(i_k+1).
\end{eqnarray*}
It is straightforward to verify that $i_k=\Theta(f(k)\log^{-1}n)$. Since $|S_k^-|=\Theta(f(k))$, the number of vertices in $S_k^-$ that \emph{do not} have degree at least $k$ is $O(|S_k^-|\log^{-1}n)$.

Next, consider the vertices in $S_k^+$. From the bounds on $k$, it follows that $R_i^*= \Omega(\omega(n)\log^2(n))$. If $R_i\leq (1-\log^{1-s}n)i$, then $R^*_i/R_i=O(\log n)$, and thus we may assume that $R_i\geq \omega(n)$.

If $R_i^*\leq n\log^{-3/\alpha}n$ and $\omega(n)\leq R_i\leq (1-\log^{1-s}n)i$ then, by Theorem \ref{thm:deg}, \emph{wep} $\deg(v_i,n)<k$.  If $R^*_i>n\log^{-3/\alpha}n$, then, by the second part of Theorem \ref{thm:deg}, \emph{wep} $\deg(v_i,n)=O(\log^{3+\alpha}n)<k$.

This time, let $f(k)=\big(1+\log^{1-s} n \big) n \left(\frac {1-\alpha}{\alpha} \cdot \frac dk \right)^{1/\alpha}$. If $R_i>(1-\log^{1-s}n)i$ and $i\geq 3f(k)$ then we may assume that $R_i\geq 2f(k)$, and, using the second part of Theorem \ref{thm:deg}, we find that \emph{wep} $\deg(v_i,n)<k$. The probability that $R_i>(1-\log^{1-s}n)i$ equals $1-(1-\log^{1-s}n)^s\leq s\log^{1-s}n$. So the expected number of vertices $v_i$ with $i\leq 3f(k)$ and $R_i>(1-\log^{1-s}n)i$ is $O(\log^{1-s}{n}f(k))$, and, using the Chernoff bound again, we can conclude that \emph{wep} the actual number is at most of the same order. Since $f(k)$ and $|S_k^-|$ are both $\Theta(k^{-1/\alpha}n)$, we can conclude that the total number of vertices in $S_k^+$ that \emph{do not} have degree less than $k$ is $O(|S_k^-|\log^{1-s}n)$. This completes the proof of the theorem.
\end{proof}

If $s<1$, then the solution of the differential equation (\ref{z:de}) is the same as for the case that $s<1$. Using methods almost identical to the ones used for the case where $s>1$, we can show that \emph{wep} the rank is close to the one suggested by the differential equation.

\begin{theorem}
For all $i\geq n^{1/2}\log^{s+1}n$, if vertex $v_i$ has initial rank $R_i$ so that $n^{1/2}\log^{s+1}{n}\leq R_i$, then \emph{wep}
\[
r(v_i,t)= \left(t - \frac{(R_i^*)^{1-s}}{1-s}t^s\right)(1+O(t^{s-1}n^{(1-s)/2}\log n)).
\]
\end{theorem}

\begin{proof}
The initial part of the proof is identical to the proof of Theorem \ref{thm:rank}, and is thus omitted. Using the differential equation method, we can show that

\begin{eqnarray*}
r(v_i,t)^{1-s}&=&(R^*_i)^{1-s}+(t+1)^{1-s}+O(R_i^{-s}n^{1/2}\log{n})\\
&=& \left( (R^*_i)^{1-s}+(t+1)^{1-s}\right)(1+O(t^{s-1} R_i^{-s}n^{1/2}\log{n}))\\
&=& \left( (R^*_i)^{1-s}+(t+1)^{1-s}\right)(1+O(t^{s-1} n^{(1-s)/2}\log{n}))
\end{eqnarray*}
for $i+1 \le t \le n$, so
\begin{eqnarray*}
r(v_i,t)&=& t\left( \left(\frac{R_i^*}{t}\right)^{1-s}+\left(\frac{t+1}{t}\right)^{1-s}\right)^{\frac{1}{1-s}}(1+O(t^{s-1} n^{(1-s)/2}\log{n}))\\
&=& \left(t+\frac{(R_i^*)^{1-s}}{1-s}t^s\right)(1+O(t^{s-1} n^{(1-s)/2}\log{n})).
\end{eqnarray*}
\end{proof}

Thus, the rank of a vertex at time $t$ tends to be close to $t$, which means we are in a situation similar to the inverse age case. In emulation of Theorem~\ref{thm:inverse_age}, we can show that almost all vertices have expected degree at most $O(\log n)$. Since the proof is almost identical to the proof of Theorem~\ref{thm:inverse_age}, it is omitted.

\begin{theorem}
Let $0<\alpha<1$, $d \in \N$, and $i \geq n^{1/2}\log^{2(1-s)}n$. Then
\begin{eqnarray*}
\E \deg(v_i,n) =O(\log n).
\end{eqnarray*}
\end{theorem}

\section{Acknowledgements}

The authors would like to thank
William Aiello and O-Yeat Chan for helpful discussions on the topics of this paper, and Filippo Menczer for suggesting the problem during WAW 2007.

\end{document}